\documentclass[
11pt,%
tightenlines,%
twoside,%
onecolumn,%
nofloats,%
nobibnotes,%
nofootinbib,%
superscriptaddress,%
noshowpacs,%
centertags]%
{revtex4}
\usepackage{ljm}
\newtheorem{proposition}{Proposition}
\begin{document}

\titlerunning{Almost quasi-Sasakian manifolds} 
\authorrunning{S.V. Galaev} 

\title{Almost quasi-Sasakian manifolds equipped\\ with skew-symmetric connection}

\author{\firstname{S.~V.}~\surname{Galaev}}
\email[E-mail: ]{sgalaev@mail.ru}
\affiliation{Saratov State University,  ul. Astrakhanskaya 83, Saratov, 410012 Russia}

\firstcollaboration{(Submitted by A.~A.~Editor-name)} 

\received{June 13, 2018} 

\begin{abstract} 
		 On a sub-Riemannian manifold, a connection with skew-symmetric torsion is defined as the unique connection from the class of $N$-connections that has this property. Two cases are considered separately: sub-Riemannian structure of even rank, and sub-Riemannian structure of odd rank. The resulting connection, called the canonical connection, is not a metric connection in the case when the sub-Riemannian structure is of even rank. The structure of an almost quasi-Sasakian manifold is defined as an almost contact metric structure of odd rank that satisfies additional requirements. Namely, it is required that the canonical connection is a metric connection and that the transversal structure is a K\"ahler structure. Both the quasi-Sasakian structure and the more general almost contact metric structure, called an almost quasi-Sasakian structure, satisfy these requirements. Sufficient conditions are found for an almost quasi-Sasakian manifold to be an Einstein manifold.
\end{abstract}



\maketitle


\section{Introduction}

 By an almost quasi-Sasakian manifold (AQS-manifold) we mean an almost normal almost contact metric manifold with a closed fundamental form for which the condition $$ d \eta   ( \vec{\xi} , \cdot  ) = 0 $$ holds true. An almost contact metric manifold is called by the author of this paper  almost normal if the equality $$ {\tilde {N}} _ {\varphi} = N _ {\varphi} +2 {\varphi} ^ * d \eta \otimes  \vec{\xi} = 0 $$ holds [5]. The ,,almost normality'' condition is equivalent to the integrability of the transversal structure. An AQS-manifold is a generalization of a quasi-Sasakian manifold and reduces to the latter if $$ d \eta = {\varphi } ^ * d \eta.$$  An AQS-structures naturally arises on the distribution of an almost contact metric manifold [3]. In this paper, an AQS-structure is equipped with a connection with skew-symmetric torsion. The study of such connections (called canonical in this work) is motivated by the demands of theoretical physics [2,4]. Initially, an $N$-connection $\nabla ^N $ is defined on a sub-Riemannian manifold of contact type endowed with an endomorphism $N: TM \to TM$ of the tangent bundle of the manifold $M$ ($N \vec{\xi} = \vec {0}$, $ N (D) \subset D $). The connection $\nabla ^N $  expressed as follows in terms of the Levi-Civita connection $\widetilde  \nabla$:
\[ \nabla ^ N_XY = {\widetilde \nabla} _ XY +   ({\widetilde \nabla} _ X \eta  )   ( Y  )  \vec{\xi} - \eta   (Y  ) {\widetilde \nabla} _ X  \vec{\xi} - \eta   (X  )   ({\widetilde \nabla} _ { \vec{\xi}} \eta  )   (Y  )  \vec{\xi} -\eta   (X  ) (C + \psi -N) Y. \]

We divide sub-Riemannian structures into two classes:  structures of even rank and structures of odd rank. For structures of odd rank, the equality $$ {\widetilde  \nabla} _ { \vec{\xi}} \,\eta = 0 $$ is satisfied, in this case the expression for the $N$-connection becomes simpler:
\[ \nabla ^ N_XY = {\widetilde \nabla} _ XY +   ({\widetilde \nabla} _ X \eta  )   ( Y  )  \vec{\xi} - \eta   (Y  ) {\widetilde \nabla} _ X  \vec{\xi} - \eta   (X  ) ( C + \psi -N) Y. \]

The case $ {\widetilde  \nabla} _ { \vec{\xi}}\, \eta = 0 $ is considered in detail in [1]. For the corresponding manifolds, conditions are found under which the $N$-connection $  \nabla ^ N $ has skew-symmetric torsion. This connection is uniquely defined and corresponds to the endomorphism $ N = 2 \psi $, where the endomorphism $ \psi $ is given by the equality $$ \omega   (X, Y  ) = g   (\psi X, Y  ) $$ and is called in the paper  the second structural endomorphism of an almost contact metric manifold. In the same paper, the concept of an $\nabla ^ N$-Einstein almost contact metric manifold is introduced. For the case $ N = 2 \psi $, conditions are found under which an almost contact metric manifold is a $\nabla ^N $-Einstein manifold.

This work consists of two main parts. In the first part, conditions are found under which an $N$-connection on a sub-Riemannian manifold of contact type has skew-symmetric torsion. In the second part,  $N$-connections with skew-symmetric torsion are considered applied to AQS-manifolds. Conditions are found under which AQS-manifolds are Einstein manifolds with respect to the canonical connection. Examples of such manifolds are given.

\section{Sub-Riemannian manifolds equipped with a canonical skew-symmetric connection}

Let $M$ be a smooth manifold of odd dimension   $ n = 2m + 1 $ \textit {} with a sub-Riemannian structure $ (M,  \vec{\xi}, \eta, g, D) $ of a contact type, where $ \eta $ and $  \vec{\xi} $ is a 1-form and a unit vector field generating, respectively, mutually orthogonal distributions $ D $ and $ D ^ {\bot} $.

Throughout the work we use adapted coordinates. A chart $k(x ^ i)$, $i, j, k = 1, \dots, n$, of the manifold $ M $ will be call adapted to the  distribution $ D $ if $$D^\bot=\textrm{Span}\left(\partial_n\right),\quad {\partial}_n  =  \vec{\xi}, $$ see [3]. 
Let $ P: TM \to  D $ be the projector defined by the decomposition $$ TM = D \oplus D ^ {\bot } ,$$ and $ k (x ^ i) $ be an adapted chart. The vector fields $$ P   ({\partial } _a  ) = \vec{e} _ a = {\partial} _a-\Gamma ^ n_a {\partial} _n,\quad a, b, c = 1, \dots, n-1, $$ are linearly independent and linearly generate the distribution  $ D = \textrm{Span} (\vec{e} _ a) $ in the definition domain of the corresponding chart.

For the adapted charts $ k (x ^ i) $ and $ k'(x ^ {i'}) $ the following coordinate transformation formulas are fulfilled:
\[x ^ a = x ^ a  (x ^ {a'}  ),\quad x ^ n = x ^ {n'} + x ^ n (x ^ {a'} ). \]

A tensor field $t$ of type $ (p, q) $ defined on an almost contact metric manifold is called admissible (to the distribution $ D $) or transversal if $t$ vanishes whenever one of its arguments is $  \vec{\xi} $ or $ \eta $. The coordinate representation of an admissible tensor field with respect to an adapted chart is as follows:
\[t = t ^ {a_1 \dots a_p} _ {b_1 \dots b_q} \vec{e} _ {a_1} \otimes \dots \otimes \vec{e} _ {a_p} \otimes dx ^ {b_1} \otimes \dots \otimes dx ^ {b_q}. \]
The transformation of the components of an admissible tensor field in adapted coordinates obeys the following law:
$$ t ^ a_b = A ^ a_ {a '} A ^ {b '} _ bt ^ {a '} _ {b '} ,$$ where $ A ^ {a '} _ a = \frac{\partial x ^ {a '}} {\partial x ^ a} $.

Thus  adapted coordinates play the role of ,,holonomic'' coordinates for the non-involutive distribution. Adapted coordinates are essentially used in foliation geometry [7]. 

Let $\omega=d\eta$. The equality $$[\vec{e} _ a,  \vec{e} _ b] = 2 {\omega } _ {ba} {\partial} _n$$ holds true. This, in particular, implies an important statement for what follows: the condition $ d \eta   ( \vec{\xi}, X  ) = 0 $ is equivalent to the equality $ {\partial} _n {\mathrm {\Gamma}} ^ n_a = 0 $.

Let $ \widetilde  \nabla $ be the Levi-Civita connection and $ {\tilde {\Gamma}} ^ i_ {jk} $ its Christoffel symbols. The following proposition may be obtained by direct computations based on the use of the equality
\[2\Gamma ^ m_ {ij} = g ^ {km}   (\vec e_ig_ {jk} + \vec e_jg_ {ik} -\vec e_kg_ {ij} + \Omega ^ l_ {kj} g_ {li} + \Omega ^ l_ {ki} g_ {lj}  ) + \Omega ^ m_ {ij}. \]

\begin{proposition} The Christoffel symbols $ {\tilde {\Gamma}} ^ k_ {ij} $ of the Levi-Civita connection of a sub-Riemannian manifold in adapted coordinates have the form:
\begin{align*}
 {\tilde {\Gamma}} ^ c_ {ab} &= \Gamma ^ c_ {ab},\quad  {\tilde {\Gamma}} ^ n_ {ab} = {\omega } _ {ba} -C_ {ab},\quad {\tilde {\Gamma}} ^ b_ {an} = {\tilde {\Gamma}} ^ b_ {na} = C ^ b_a + {\psi} ^ b_a,\\  {\tilde {\Gamma}} ^ n_ {na} &= - {\partial} _n\Gamma ^ n_a,\quad  {\tilde {\Gamma}} ^ a_ {nn} = g ^ {ab} {\partial} _n\Gamma ^ n_b,\end{align*} where $$\Gamma ^ a_ {bc } = \frac {1} {2} g ^ {ad} (\vec{e} _ bg_ {cd} + \vec{e} _ cg_ {bd} - \vec{e} _ dg_ {bc }),\quad  {\psi } ^ b_a = g ^ {bc} {\omega} _ {ac},\quad  C_ {ab} = \frac {1} {2} {\partial} _ng_ {ab}, \quad C ^ b_a = g ^ {bc} C_ {ac}.$$ Here the endomorphism $ \psi: TM \to TM $ is determined by the equality $$ \omega   (X, Y  ) = g   (\psi X, Y  ) ,$$ and we set $$ C   (X, Y  ) = \frac {1} {2} (L _ { \vec{\xi} }\, g)   (X, Y  ). $$
\end{proposition}

An $N$-connection $  \nabla ^ N $ is defined on a sub-Riemannian manifold endowed with the endomorphism $ N: TM \to TM $ of the tangent bundle of $ M $ ($ N  \vec{\xi} = \overrightarrow {0 } $, $ N (D) \subset D $). The connection $\nabla^N$ may be  expressed in terms of the Levi-Civita connection $ {\widetilde  \nabla}$,
\[ \nabla ^ N_XY = {\widetilde \nabla} _ XY +   ({\widetilde \nabla} _ X \eta  )   ( Y  )  \vec{\xi} - \eta   (Y  ) {\widetilde \nabla} _ X  \vec{\xi} - \eta   (X  )   ({\widetilde \nabla} _ { \vec{\xi}} \eta  )   (Y  )  \vec{\xi} - \eta   (X  ) (C + \psi -N) Y. \]

\begin{proposition}  A linear connection $  \nabla ^ N $ defined on a sub-Riemannian manifold is skew-symmetric if and only if $ N = 2 \psi $.
\end{proposition}

{\bf Proof.} It may be  directly checked that with respect to  adapted coordinates the nonzero Christoffel symbols $ G ^ i_ {jk} $ of the connection $  \nabla ^ N _ {\vec{x}} $ have the form $$ G ^ a_ {bc} = \frac {1} {2} g ^ {ad}   (\vec{e} _ bg_ {cd} + \vec{e} _ cg_ {bd} - \vec{e} _ dg_ {bc}  ),\quad  {\mathrm {G}} ^ b_ {na} = N ^ b_a, \quad {\mathrm {G}} ^ n_ {na}  = - {\partial} _n\Gamma ^ n_a .$$

The rank of a sub-Riemannian structure is  equal to $ 2p $ if $ {(d \eta ) } ^ p \neq 0 $, $ \eta \wedge {(d \eta)} ^ p = 0 $, and equal to $ 2p + 1 $ if $ \eta \wedge {(d \eta)} ^ p \neq 0 $, $ {(d \eta)} ^ {p + 1} = 0 $. It is easy to check that the rank of a sub-Riemannian structure is $ 2p + 1 $ if and only if $ {\partial} _n\Gamma ^ n_a = 0 $.

Put $ \tilde {S}   (X, Y, Z  ) = g (S   (X, Y  ), Z) $, $ X, Y, Z \mathrm {\in } \mathrm {\Gamma} \mathrm {(} \mathrm {TM} \mathrm {)} $. With respect to  adapted coordinates the nonzero components of the tensor $ \tilde {S}   (X, Y, Z  ) $  have the following form:
\[\tilde {S}   (\vec{e} _ a, \vec{e} _ b, {\partial } _n  ) = 2 {\omega} _ {ab}, \]
\[\tilde {S}   (\vec{e} _ a, {\partial } _n, \vec{e} _ b  ) = - g (N \vec{e} _ a, \vec{e} _ b), \]
\[\tilde {S}   ({\partial } _n, \vec{e} _ a, \vec{e} _ b  ) = g (N \vec{e} _ a, { \vec{e}} _ b). \]

The tensor $ \tilde {S}   (X, Y, Z  ) $ is skew-symmetric if and only if $ 2 {\omega } _ {ab} = g (N \vec{e} _ a, \vec{e} _ b) $. This proves Proposition 2. \qed

If $ N = 2 \psi $, the connection $  \nabla ^ N $ will be called the canonical connection.
Note that the canonical connection in the case of a sub-Riemannian structure of even rank is not a metric connection. Indeed, $$  \nabla ^ N_ng_ {na} = - {\mathrm {G}} ^ n_ {na} = {\partial} _n\Gamma ^ n_a .$$

\section{Basic information from the geometry of almost quasi-Sasakian manifolds}

Consider an almost contact metric manifold $ M $ of odd dimension   $ n = 2m + 1 $. Let $(M,  \vec{\xi}, \eta, \varphi, g, D  ) $ be an almost contact metric structure on a manifold $ M $, where $ \varphi $ is a tensor of type~$(1 , 1)$, called a structural endomorphism, $  \vec{\xi} $ and $ \eta $ are a vector and a covector, called, respectively, a structure vector and a contact form, $g$ is a (pseudo-)Riemannian metric. In this case, the following equalities hold true:

\begin{itemize}
\item[1)] $ {\varphi } ^ 2 = -I + \eta \otimes  \vec{\xi} $,
\item[2)] $ \eta   ( \vec{\xi}  ) = 1 $,
\item[3)] $ g   (\varphi X, \varphi Y  ) = g   (X, Y  ) - \eta   (X  ) \eta (Y) $,  $ X, Y \in \mathit {\Gamma} (TM) $. 
\end{itemize}

The smooth distribution $ D = \mathrm {ker}  (\eta ) $ is called the distribution of an almost contact structure.

As a consequence of conditions 1) - 3) we obtain:
\begin{itemize}
\item[5)] $\varphi  \vec{\xi} = \vec {0},$
\item[6)] $\eta \circ \varphi = 0,$
\item[7)] $\eta   (X  ) = g (X,  \vec{\xi}),$ $X \in {\Gamma} (TM)$. \end{itemize}

The skew-symmetric tensor $ \mathrm {\Omega }   (X, Y  ) = g (X, \varphi Y) $ is called the fundamental form of the structure. An almost contact metric structure is called a contact metric structure if the equality $ \mathrm {\Omega } = d \eta $ holds. The smooth distribution $ D ^ {\bot} = \mathrm{Span} ( \vec{\xi}) $, orthogonal to the distribution $ D $, is called the framing of the distribution $ D $. There is the decomposition $ TM = D \oplus D ^ {\bot} $.

A Sasakian manifold is a contact metric space satisfying the additional condition $$ N ^ {(1)} _ {\varphi } = N _ {\varphi} + 2d \eta \otimes  \vec{\xi} = 0,$$ where $$ N _ {\varphi}   (X, Y  ) =   [\varphi X, \varphi Y  ] + {\varphi} ^ 2   [X, Y  ] - \varphi   [\varphi X, Y  ] - \varphi [X, \varphi Y] $$
	 is the Nijenhuis tensor of the endomorphism $\varphi$. The   condition $ N ^ {(1)} _ \varphi  = N _ \varphi + 2d \eta \otimes  \vec{\xi} = 0 $ means that the Sasakian space is a normal space.

An almost contact metric manifold is called an almost contact K\"ahler manifold [5] if the following conditions hold: $$ d\Omega = 0 ,\quad  {\tilde {N}} _ {\varphi } = N _ {\varphi} +2 {\varphi} ^ * d \eta \otimes  \vec{\xi} = 0 .$$ A manifold for which the condition $ {{\tilde {N}} _ {\varphi} = N} _ {\varphi} +2 {\varphi} ^ * d \eta \otimes  \vec{\xi} = 0 $ is satisfied is called by us an almost normal manifold. It is easy to check that an almost normal almost contact metric manifold is a normal manifold if and only if $ d \eta = {\varphi } ^ * d \eta $.

Let $ P: TM \to D $ be the projector defined by the decomposition $ TM = D \oplus D ^ {\bot } $. Then the following proposition holds.

\begin{proposition} For any almost contact metric manifold, the following equality holds: $ PN ^ {(1)} _ {\varphi } = {\tilde {N}} _ {\varphi} $.\end{proposition}

{\bf Proof.} \begin{multline*} PN ^ {(1)} _ {\varphi }   (X, Y  ) = P   (  [\varphi X, \varphi Y  ] + {\varphi} ^ 2   [X, Y  ] - \varphi   [\varphi X, Y  ] - \varphi   [X, \varphi \vec {y}  ] + 2d \eta   (X, Y  )  \vec{\xi}  )\\ = P   [\varphi X, \varphi Y  ] -P   [X, Y  ] - \varphi   [\varphi X, Y  ] - \varphi   [X, \varphi Y  ]\\ =   [\varphi X, \varphi Y  ] - \eta   (  [\varphi X, \varphi Y  ]  )  \vec{\xi} -   [X, Y  ] + \eta   (  [X, Y  ]  )  \vec{\xi} - \varphi   [\varphi X, Y  ] - \varphi   [X, \varphi Y  ] \\ =   [\varphi X, \varphi Y  ] + {\varphi} ^ 2   [X, Y  ] - \varphi   [\varphi X, Y  ] - \varphi   [X, \varphi Y  ] + 2d \eta   (\varphi X, \varphi Y  )  \vec{\xi} = {\tilde {N}} _ {\varphi}   (X, Y  ) . \end{multline*}

The proposition is proved. \qed 

Note that the  just proved proposition implies the relation:
\begin {equation} \label {GrindEQ__1_}
N ^ {(1)} _ {\varphi }   (X, Y  ) = {\tilde {N}} _ {\varphi}   (X, Y  ) +2 (d \eta   (X, Y  ) -d \eta   (\varphi X, \varphi Y  ))  \vec{\xi}.
\end {equation}

An internal linear connection $ \mathrm {\nabla } $ [3] on a manifold with an almost contact metric structure is a map
\[\mathrm {\nabla } : \mathrm {\Gamma} \mathrm {(} D) \times \mathrm {\Gamma} \mathrm {(} D) \to \mathrm {\Gamma} \mathrm {(} D), \]

satisfying the following conditions:
\begin{itemize}
	\item[1)]  $\nabla _ {f_1 \vec {x} + f_2 \vec {y}} = f_1 \nabla _ {\vec{x}} + f_2 \nabla _ {\vec {y}}, $
\item[2)]  $\nabla _ {\vec{x}} f \vec {y} =   (\vec {x} f  ) \vec {y} + f {\mathrm { \nabla}} _ {\vec{x}} \vec {y},$
\item[3]  $\nabla _ {\vec{x}}   (\vec {y} + \vec {z}  ) = \nabla _ {\vec{x}} \vec {y} + \nabla _ {\vec{x}} \vec {z},$ \end{itemize}

where ${\Gamma } (D)$ is the module of admissible vector fields (vector fields at each point belonging to the distribution $ D $).

The Christoffel symbols of $\nabla$ are determined from the relation $  \nabla _ {\vec{e} _ a} \vec{e} _ b = \Gamma ^ c_ {ab} { \vec {e}} _ c $. From the equality $ \vec{e} _ a = A ^ {a '} _ a \vec{e} _ {a '} $, where $ A ^ {a '} _ a = \frac {\partial x ^ {a '}} {\partial x ^ a} $, the transformation formula  follows
\[{\Gamma} ^ c_ {ab} = A ^ {a '} _ aA ^ {b '} _ bA ^ c_ {c '} \Gamma ^ {c '} _ {a ' b '} + A ^ c_ {c '} \vec{e} _ aA ^ {c '} _ b.\]
Hence, in particular, it follows that the derivatives $ {\partial } _n {\mathrm {\Gamma}} ^ d_ {ac} $ are components of an admissible tensor field.

We give two  examples of almost contact K\"ahler manifolds.

\textbf{Example 1.} Let $ M = \{(x, y, z, u, v) \in R ^ 5: y \neq 0 \} $ be a smooth manifold of dimension 5 equipped with an almost contact
metric structure $  (M,  \vec{\xi} , \eta, \varphi, g, D  ) $, where
\begin{itemize}
\item[1)] $ D = \mathrm{Span}\{\vec{e} _ {1}, \vec{e} _ {2}, \vec{e} _ {3}, \vec{e} _ {4}\} $, here $ \vec{e} _ 1 = {\partial} _1-y {\partial} _5$, $\vec{e} _ 2 = {\partial} _2,$ $\vec{e} _ 3 = {\partial } _3,$  $\vec{e} _ 4 = {\partial} _4 $, and $(\partial_1,\dots,\partial_5)$ is the basis of vector fields corresponding to the coordinates $(x,y,z,u,v)$ on $\mathbb{R}^5$, 
\item[2)] $  \vec{\xi} = {\partial} _5 $, 
\item[3)] $ \eta = dz + ydx $, 
\item[4)] $ \varphi \vec{e} _ 1 = \vec{e} _ 3,$  $\varphi \vec{e} _ 2 = {\vec {e }} _ 4,$ $\varphi \vec{e} _ 3 = - \vec{e} _ 1,$ $\varphi \vec{e} _ 4 = - \vec{e} _ 2,$ $\varphi  \vec{\xi} = 0$, \item[5)] the basis $ (\vec{e} _ {1}, \vec{e} _ {2}, \vec{e} _ {3}, \vec{e} _ {4,  }  \vec{\xi}) $
consists of orthonormal vectors. 
\end{itemize}
It may be directly checked that the almost contact metric manifold $ M $ is not normal, but almost normal. Indeed, $$ N ^ {  (1  )} _ {\varphi }   (\vec{e} _ 1, \vec{e} _ 2  ) = {\varphi} ^ 2   [\vec{e} _ 1, \vec{e} _ 2  ] +  [\vec{e} _ 3,  \vec{e} _ 4  ] - \varphi   [\vec{e} _ 3,  \vec{e} _ 2  ] - \varphi   [\vec{e} _ 1,  \vec{e} _ 4  ] + 2d \eta   (\vec{e} _ 1,  \vec{e} _ 2  )  \vec{\xi} = {\varphi} ^ 2  \vec{\xi} - \eta  ( \vec{\xi}  )  \vec{\xi} = -  \vec{\xi} .$$ 
On the other hand, $$ {\tilde {N}} _ {\varphi}   (\vec{e} _ 1, \ \vec{e} _ 2  ) = 2d \eta   (\vec{e} _ 3, \ \vec{e} _ 4  )  \vec{\xi} = 0 .$$
For the structure under consideration, the equality $ d \eta   ( \vec{\xi}, X  ) = 0 $, $ X \in \ {\Gamma} (TM) $, holds. 
Thus $ \omega = d \eta $ is an admissible tensor field, to which the internal connection $\nabla $ may be applied [3]. 
Moreover, $  {\nabla } \omega = 0 $. Let, further, $ \psi $ be the endomorphism defined by the equality $ \omega   (X, Y  ) = g   (\psi X, Y  ) $. The coordinate representation of the endomorphism $ \psi$ is of the form $ {\psi} ^ b_a = g ^ {bc} {\omega} _ {ac} $.
Thus  the trace of the square of the endomorphism $ \psi $ is constant, $ \mathrm{tr}   ({\psi} ^ 2  ) = const $. 

\textbf {Example 2. } Considers the same manifold $ M $ as in the previous example with the only difference that $$ \vec{e} _ 1 = {\partial } _1-y {z \partial} _5 , \quad \eta = dz + yzdx. $$ Unlike the previous case, the condition $ d \eta   ( \vec{\xi}, \cdot  ) = 0 $ is not satisfied. Indeed, $$ 2d \eta   ( \vec{\xi} , \vec{e} _ 1  ) = - \eta (  [ \vec{\xi}, \vec{e} _ 1  ]) = y \neq 0 .$$

An almost contact metric manifold is called an almost quasi-Sasakian manifold (AQS-manifold) if the following conditions are satisfied: \begin{equation}\label{condAQS} d\Omega = 0 , \quad  {\tilde {N}} _ {\varphi } = N _ {\varphi} +2 {\varphi} ^ * d \eta \otimes  \vec{\xi} = 0 , \quad  d \eta   ( \vec{\xi}, \cdot  ) = 0. \end{equation} 

Note that Example 1 implies that there exist an almost quasi-Sasakian manifolds  satisfying the conditions $$ \nabla \omega = 0 ,\quad  \mathrm{tr}   ({\psi} ^ 2  ) = const .$$ Moreover the equality $ \mathrm {\nabla} \omega = 0 $ is equivalent to the equality $  {\nabla} \psi = 0. $

The following theorem holds.

\begin{theorem} An almost contact metric structure is an almost quasi-Sasakian structure if and only if the following equality holds:
\begin {equation} \label {GrindEQ__2_}
{(\widetilde  \nabla} _ X \varphi) Y = g   (  (\psi \circ \varphi  ) Y, X  )  \vec{\xi} - \eta   (Y  )   (\varphi \circ \psi  )   (X  ) - \eta   (X  )   (\varphi \circ \psi - \psi \circ \varphi  ) Y.
\end {equation}\end{theorem}

{\bf Proof.}  Let $ M $ be a AQS-manifold. Let us show that the condition \eqref {GrindEQ__2_} is satisfied. The equality  
\begin{multline*}2g   ({(\widetilde  \nabla} _ X \varphi  ) Y, Z) = 3   (d\Omega   (\mathrm {X,} \ \varphi Y, \varphi Z  ) -d\Omega   (\mathrm {X,} \ Y, Z  )  ) + g   (N ^ {  (1  )} _ {\varphi}   ( Y, Z  ), \varphi X  ) \\+ {2N} ^ {  (2  )} _ {\varphi}   (Y, Z  ) \eta   (X  ) +2 (d \eta   (\varphi Y, X  ) \eta   (Z  ) -d \eta   (\varphi Z, X  ) \eta   (Y  ) ), \end{multline*} where
the operator $N^{(2)}_\varphi$ is similar to $N^{(1)}_\varphi$, see, e.g., [1].
holds for any almost contact metric manifold.
Using \eqref{condAQS}, we get $$ g   ({(\widetilde \nabla} _ X \varphi  ) Y, Z) = \eta   (X  )   (g   (  (\psi \circ \varphi  ) Y, Z  ) + d \eta   (Y, \varphi Z  )  ) + g   (d \eta   (\varphi Y, X  )  \vec{\xi} , Z  ) + d \eta   (X, \varphi Z  ) \eta   (Y  ).$$ This proves the equality \eqref {GrindEQ__2_}. The inverse statement may be easily proved  
using adapted coordinates. \qed

The following propositions are direct consequences of Theorem 1.

\begin{proposition} An almost contact metric structure is a quasi-Sasakian structure if and only if the equality holds:
\[{(\widetilde \nabla} _ X \varphi) Y = g   (AY, X  )  \vec{\xi} - \eta   (Y  ) AX, A = \varphi \circ \psi. \]
\end{proposition}

\begin{proposition} An almost quasi-Sasakian manifold is a quasi-Sasakian manifold if and only if one of the following conditions holds
	\begin{itemize}
		\item[1)] $d \eta = {\varphi} ^ * d \eta,$
\item[2)] $\varphi \circ \psi - \psi \circ \varphi = 0, $
\item[3)] $g   (X, AY  ) = g (AX, Y), A = \varphi \circ \psi.$
\end{itemize}
\end{proposition}

\section{Almost quasi-Sasakian manifolds with canonical skew-symmetric connection}

Let $   (M,  \vec{\xi} , \eta, \varphi, g, D  ) $ be an almost quasi-Sasakian structure given on a manifold $ M $, and let $ \nabla ^ N $ be the canonical connection. Proposition 5 implies the following 

\begin{proposition} An almost quasi-Sasakian manifold is a quasi-Sasakian manifold if and only if $  \nabla ^ N \varphi = 0 $. \end{proposition}

Further, we restrict our attention to the case when the torsion $ \tilde {S}   (X, Y, Z  ) $ of the connection $  \nabla ^ N $ is parallel. It is known that this condition holds for  Sasakian manifolds. At the same time, it follows from the above example that there exist such almost quasi-Sasakian manifolds that are neither quasi-Sasakian, nor Sasakian  and for which the   torsion is skew-symmetric and parallel.

Let $ K    $ be the curvature tensor of the canonical connection $ {\nabla}^ N $. For nonzero components of the tensor $  K  $ it holds
\[K ^ d_ {abc} = R ^ d_ {abc} +4 {\omega } _ {ab} {\psi} ^ d_c. \]
\[K ^ d_ {anc} = {\mathrm {2} \mathrm {\nabla } } _a {\psi} ^ d_c. \]

Here $ R ^ d_ {abc} = 2 \vec{e} _ {[a} {\mathrm {\Gamma } } ^ d_ {b] c} +2 {\mathrm {\Gamma}} ^ d_ { [a   \lceil e   \rceil} {\mathrm {\Gamma}} ^ e_ {b] c} $ are components of the Schouten curvature tensor [3] defined by the equality
\[R   (X, Y  ) Z =  \nabla _X \nabla _ Y Z - \nabla _ Y \nabla _ XZ - \nabla _ {P [X, Y]} Z-P [Q   [X, Y  ], Z],\quad Q = 1-P. \]
The tensor $ K $ may be written in the following way:
\[K   (X, Y  ) Z = R   (X, Y  ) Z + \eta   (Y  )   ( \nabla _XN  ) Z - \eta   (X  )   (\nabla _ YN  ) Z + 4 \omega   (X, Y  ) \psi (Z),\quad X, Y, Z \mathrm {\in} \Gamma \mathrm {(TM)}. \]

Tensor field $ r   (X, Z  ) =\mathrm {tr} ( Y \mapsto  \mathrm {R}   (X, Y  ) Z) $, $ X, Z \mathrm {\in} \Gamma \mathrm {(D)} $, will be called the Ricci-Wagner tensor.
Using  adapted coordinates, we write down the components of the Ricci tensor $ k  $ of the connection~$  \nabla ^ N $:
\[k_ {ab} = r_ {ab} +4 {\omega } _ {ad} {\psi} ^ d_b, \]
\[k_ {an} = k_ {nn} = 0, \]
\[k_ {na} = -  \nabla _d {\psi} ^ d_a. \]

From the definition of the  endomorphism $ \psi $ it follows that the equality $ \mathrm {\nabla } \omega = 0 $ holds if and only if $ {\nabla} \psi = 0 $. Hence the components of the Ricci tensor $ k$ of the connection $ \nabla ^ N $ with parallel torsion take the form:
\[k_ {ab} = r_ {ab} +4 {\omega } _ {ad} {\psi} ^ d_b, \]
\[k_ {an} = k_ {na} = k_ {nn} = 0. \]

Thus the following theorem turns out to be true.

\begin{theorem} An almost quasi-Sasakian manifold is an Einstein manifold with respect to the canonical connection with parallel torsion if and only if
\[r_ {ab} = 4 {\omega } _ {da} {\psi} ^ d_b. \] \end{theorem}

As follows from Theorem 2, the existence of Einstein metric on an almost quasi-Sasakian manifold with respect to the canonical connection essentially depends on the structure of the Ricci-Wagner tensor.

We complete the work with an example of an almost quasi-Sasakian Einstein manifold.

\textbf {Example 3. } Let us introduce a quasi-Sasakian structure on the manifold $ M = \{(x, y, z, u, v) \in R ^ 5: y \neq 0 \} $ by
setting: \begin{itemize}
	\item[1)] $ D = \mathrm{Span}\{\vec{e} _ {1}, \vec{e} _ {2}, \vec{e} _ {3}, \vec{e} _ { 4 \ }\}$, 				
where $ \vec{e} _ 1 = {\partial} _1-y {\partial} _5, \  \vec{e} _ 2 = {\partial} _2, \ \vec{e} _ 3 = {\partial} _3, \ \vec{e} _ 4 = {\partial} _4 $,  and $(\partial_1,\dots,\partial_5)$ is the basis of vector fields corresponding to the coordinates $(x,y,z,u,v)$ on $\mathbb{R}^5$,
\item[2)] $  \vec{\xi} = {\partial} _5 $, \item[3)] $ \eta = dv + ydx $, \item[4)] $ \varphi \vec{e} _ 1 = \vec{e} _ 2, \ \varphi \vec{e} _ 2 = - \vec{e} _ 1, \ \varphi \vec{e} _ 3 = \vec{e} _ 4,  \varphi \vec{e} _ 4 = - \vec{e} _ 3, \ \varphi  \vec{\xi} = 0 $,
\item[5)] the
metric tensor is given by the equality $$ g = \frac {1} {{  (1 + x ^ 2 + y ^ 2  )} ^ 2}   ({  (dx  )} ^ 2 + {  (dy  )} ^ 2  ) + {  (dz  )} ^ 2 + {  (du  )} ^ 2 + {  \eta} ^ 2 .$$
\end{itemize} The condition $ r_ {ab} = 4 {\omega} _ {da} {\psi} ^ d_b $  reduces to the equality $ {r _ {\alpha \beta } = -4g} _ {\alpha \beta} $,  $ \alpha, \beta = 1,2 $. It is easy to check that the above metric satisfies the condition $ {r _ {\alpha \beta} = - 4g} _ {\alpha \beta} $. Thus, we have obtained an example of an Einstein quasi-Sasakian manifold. If we redefine the first structural endomorphism by setting $$ \varphi \vec{e} _ 1 = \vec{e} _ 3, \ \varphi \vec{e} _ 2 = \vec{e} _ 4, \ \varphi \vec{e} _ 3 = - \vec{e} _ 1, \varphi \vec{e} _ 4 = - \vec{e} _ 2, \varphi  \vec{\xi} = 0 ,$$ then the quasi-Sasakian manifold reduces to an almost quasi-Sasakian manifold. Thus we have obtained an example of an Einstein AQS-manifold with respect to the canonical connection.


\begin{thebibliography}{99}
	
	\bibitem{1} S.V. Galaev,  $  \nabla ^ N $-Einstein almost contact metric manifolds. Bulletin of Tomsk State University. Mathematics and Mechanics 2021, no. 70,  5-15.

\bibitem{2}  I.Agricola, A.C. Ferreira, Einstein manifolds with skew torsion, Qaur. J. Math. 65 (2014), no. 3, 717-741.

\bibitem{3} A.V. Bukusheva, S.V. Galaev, Almost contact metric structures defined by connection over distribution. Bulletin of the Transilvania University of Brasov Series III: Mathematics, Informatics, Physics. 2011. Vol. 4 (2), no. 2,   13-22.

\bibitem{4} T. Friedrich,S. Ivanov, Parallel spinors and connections with skew-symmetric torsion in string theory, Asian J. Math. 6 (2002), 303--336.

\bibitem{5} S.V. Galaev, Admissible Hyper-Complex Pseudo-Hermitian Structures. Lobachevskii Journal of Mathematics, 2018, Vol. 39, no. 1, 71-76.

\bibitem{6} S.V. Galaev, Intrinsic geometry of almost contact K\"{a}hlerian manifolds.  Acta Mathematica Academiae Paedagogicae Nyiregyhaziensis. 2015. Vol 31, 35-46.

\bibitem{7} M.A. Malakhaltsev, Foliations with leaf structures, J. Math. Sci. (New York), \textbf {108}: 2 (2002), 188-210.


\end{thebibliography}
\end{document}